\documentclass[11pt]{amsart}
\usepackage{graphicx}
\usepackage{amssymb,amscd}
\newtheorem{theorem}{Theorem}[section]
\newtheorem{thm}{Theorem}[section]

\newtheorem{lemma}[theorem]{Lemma}
\newtheorem{corollary}[theorem]{Corollary}

\newtheorem{defn}[theorem]{Definition}
\newcommand{\PP}{{\mathcal P}}

\title{Convexity properties of Thompson's group $F$}

\author[Matthew Horak]{Matthew Horak}
\address{Department of Mathematics, Statistics and Computer Science, University of Wisconsin-Stout,  Menomonie, WI 54751}
\email{horakm@uwstout.edu}

\author[Melanie Stein] {Melanie Stein}
\address{Department of Mathematics, Trinity College, Hartford, CT 06106}
\email{melanie.stein@trincoll.edu}

\author[Jennifer Taback] {Jennifer Taback}
\address{Department of Mathematics, Bowdoin College, Brunswick, ME 04011}
\email{jtaback@bowdoin.edu}
\thanks{The third author acknowledges partial support from
NSF grant DMS-0604645, and the second and third authors acknowledge
partial support from a Bowdoin College Faculty Research Grant.}

\begin{document}

\maketitle

\begin{abstract}
We prove that Thompson's group $F$ is not minimally almost convex
with respect to any generating set which is a subset of the standard
infinite generating set for $F$ and which contains $x_1$.  We use
this to show that $F$ is not almost convex with respect to any
generating set which is a subset of the standard infinite generating
set, generalizing results in \cite{HST}.
\end{abstract}

\section{Introduction}

Convexity properties of a group $G$ with respect to a finite
generating set $S$ yield information about the configuration of
spheres within the Cayley graph $\Gamma(G,S)$ of $G$ with respect to
$S$. A finitely generated group $G$ is {\em almost convex(k)}, or
$AC(k)$ with respect to a finite generating set $X$ if there is a
constant $L(k)$ satisfying the following property.  For every
positive integer $n$, any two elements $x$ and $y$ in the ball
$B(n)$ of radius $n$ with $d_X(x,y) \leq k$ can be connected by a
path of length $L(k)$ which lies completely within this ball. J.
Cannon, who introduced this property in \cite{C}, proved that if a
group $G$ is $AC(2)$ with respect to a generating set $X$ then it is
also $AC(k)$ for all $k \geq 2$ with respect to that generating set.
Thus if $(G,X)$ is $AC(2)$, it is called {\em almost convex} with
respect to that generating set.

Almost convexity is a property which depends on generating set; this
was proven by C. Thiel using the generalized Heisenberg groups
\cite{T}. If a group is almost convex with respect to {\em any}
generating set, then we simply call it almost convex, omitting the
mention of a generating set. Groups which are almost convex with
respect to any generating set include hyperbolic groups \cite{C} and
fundamental groups of closed 3-manifolds whose geometry is not
modeled on {\em Sol} \cite{SS}. Moreover, amalgamated products of
almost convex groups retain this property \cite{C}.

If $(G,X)$ is not almost convex then there is a sequence of points
$\{x_i,y_i\}$ at distance $2$ in $B(n_i)$ which require successively
longer paths within $B(n_i)$ to connect them, as $i$ and $n_i$
increase. Such groups include include fundamental groups of closed
3-manifolds whose geometry is modeled on {\em Sol} \cite{CFGT} and
the solvable Baumslag-Solitar groups $BS(1,n)$ \cite{MS}, in both
cases with respect to any finite generating set, and Thompson's
group $F$ with respect to any generating set of the form $\{x_0,
x_1, \ldots, x_n\}$ which is a subset of the standard infinite
generating set for $F$ \cite{CT1,HST}.

Clearly, any two points in $B(n)$ can always be connected by a path
of length $2n$. A weaker convexity condition is {\em minimal almost
convexity}, which asks whether any two points in $B(n)$ at distance
two can be connected by a path of length at most $2n-1$ lying within
this ball. A group $G$ is said to be minimally almost convex with
respect to a finite generating set $X$ if the Cayley graph
$\Gamma(G,X)$ has this property.  In groups which are not minimally
almost convex, we can find examples of points $x,y \in B(n)$ at
distance two so that any path connecting $x$ to $y$ within $B(n)$
has length at least $2n$, even paths which do not pass through the
identity.  If $G$ is not minimally almost convex with respect to a
finite generating set $X$, then $\Gamma(G,X)$ contains isometrically
embedded loops of arbitrarily large circumference.  I. Kapovich
proved in \cite{K} that any group which is minimally almost convex
is also finitely presented.

M. Elder and S. Hermiller prove in \cite{EH} that the solvable
Baumslag-Solitar group $BS(1,2) = \langle a,t | tat^{-1} = a^2
\rangle$ is minimally almost convex with respect to the given
generating set, but for $q \geq 7$ the group $BS(1,q) = \langle a,t
| tat^{-1} = a^q \rangle$ is not minimally almost convex with
respect to the analogous generating set.  In addition, they prove
that Stallings' group:
$$S = \langle a,b,c,d,s | [a,c]=[a,d]=[b,c]=[b,d]=1, \ (a^{-1}b)^s =
a^{-1}b, \ (a^{-1}c)^s = a^{-1}c, \ (a^{-1}d)^s = a^{-1}d \rangle$$
is not minimally almost convex with respect to the above generating
set.  J. Belk and K.-U. Bux prove in \cite{BBu} that Thompson's
group $F$ is not minimally almost convex with respect to the
standard finite generating set $\{x_0,x_1\}$.

J. Meier posed a conjecture relating these two notions of convexity.
Namely, he conjectured that if a finitely generated group $G$ is not
minimally almost convex with respect to one finite generating set,
then it cannot be almost convex with respect to any finite
generating set.  We prove the following special case of this
conjecture.  Suppose $X$ and $Y$ are two finite generating sets for
a group $G$. Then $G$ can be viewed as a metric space using the
wordlength metric with respect to either generating set; we write
$(G,X)$ for $G$ viewed as a metric space using length with respect
to $X$. The identity map on $G$ is a quasi-isometry between $(G,X)$
and $(G,Y)$. We prove this conjecture in the case that this
quasi-isometry is a coarse isometry, that is, has multiplicative
constant equal to one, in Theorem \ref{thm:notmac} below.

{\bf Theorem \ref{thm:notmac}} Let $f:(G,X_G) \rightarrow (H,X_H)$ be a
$C$-coarse-isometry.
 If $(G,X_G)$ is not minimally almost convex, then $(H,X_H)$ is not almost convex.

Convexity properties have been studied for Thompson's group $F$ with
respect to its standard finite generating set $X_1=\{x_0,x_1\}$.
This group can be viewed either as a finitely or infinitely
presented group, using the two standard presentations:
$$\langle x_k, \ k \geq 0 | x_i^{-1}x_jx_i = x_{j+1} \ \text{ if
}i<j \rangle$$ or, as it is clear that $x_0$ and $x_1$ are
sufficient to generate the entire group,
 since powers of $x_0$ conjugate $x_1$ to $x_i$ for $i \geq 2$,
$$\langle x_0,x_1 |
[x_0x_1^{-1},x_0^{-1}x_1x_0],[x_0x_1^{-1},x_0^{-2}x_1x_0^2]
\rangle.$$

As noted above, the group $F$ is shown to be not almost convex with
respect to $X_1$ in \cite{CT1} and not minimally almost convex with
respect to $X_1$ in \cite{BBu}.  The proofs of these facts rely,
repectively, on the methods of computing word length in $F$ with
respect to $X_1$ due to Fordham \cite{F} and Belk and Brown
\cite{BBr}.  In \cite{HST},  we present a method for computing word
length in $F$ with respect to {\em consecutive} generating sets of
the form $X_n=\{x_0, x_1, \ldots, x_n\}$, each a finite subset of
the standard infinite generating set for $F$. This method is then
used to show that $F$ is not almost convex with respect to the
consecutive generating sets $X_n$ (Theorem 6.2, \cite{HST}).

In this paper, we first extend the result of Belk and Bux to
consecutive generating sets. We prove:

{\bf Theorem \ref{thm:notMACconsec}} Let $X_n = \{x_0,x_1, \cdots ,x_n\}$ be a consecutive generating set
for $F$.  Then $F$ is not minimally almost convex with respect to
$X_n$.

The group $F$ can be generated by any subset of the standard
infinite generating set containing $x_0$. While there are many other
finite generating sets for $F$, such as $\{x_0,x_1x_0^{-1}\}$, there
is no known method for recognizing other generating sets for this
group, or computing word length with respect to these generating
sets. We extend our initial result to show:

{\bf Theorem \ref{thm:notMACwithx_1}} Let $X=\{x_0,x_1,
x_{i_1},x_{i_2}, \cdots ,x_{i_j}\}$, where $1< i_1 < \cdots < i_j$,
be a generating set for $F$.  Then $F$ is not minimally almost
convex with respect to $X$.

We then apply Theorem \ref{thm:notmac}, the special case of J.
Meier's conjecture,  to prove:

{\bf Theorem \ref{thm:notacgeneral}} Let $X$ be any subset of the standard infinite generating set for
$F$ which includes $x_0$.  Then $F$ is not almost convex with
respect to $X$.

\section{Computing word length in Thompson's group $F$}\label{sec:computinginF}

In this section we summarize the method for computing word length of
elements of $F$ with respect to the consecutive generating sets
$X_n=\{x_0, x_1, \ldots, x_n \}$ which was introduced in \cite{HST},
and refer the reader to that paper for complete details.

Elements of $F$ can be viewed combinatorially as pairs of finite
binary rooted trees, each with the same number of carets, called
tree pair diagrams. We define a {\em caret} to be a vertex of the
tree together with two downward oriented edges, which we refer to as
the left and right edges of the caret. The {\em right (respectively
left) child} of a caret $c$ in a tree $T$ is defined to be a caret
which is attached to the right (resp. left) edge of $c$.  If a caret
$c$ does not have a right (resp. left) child, we call the right
(resp. left) edge, or leaf, of $c$ {\em exposed}.  Define the {\em
level} of a caret inductively as follows. The root caret is defined
to be at level 1, and the child of a level $k$ caret has level
$k+1$, for $k \geq 1$.

We number the leaves of each tree from $0$ through $n$, going from
left to right, and number the carets in infix order from $1$ through
$n$. The infix ordering is carried out by numbering the left child
of a caret $c$ before numbering $c$, and the right child of $c$
afterwards. Each element $g \in F$ can be represented by an
equivalence class of tree pair diagrams, among which there is a
unique reduced tree pair diagram. We say that a pair of trees is
{\em unreduced} if when the leaves are numbered from $0$ through
$n$, there is a caret in both trees with two exposed leaves bearing
the same leaf numbers.  We remove such pairs of carets, renumber the
leaves and check this condition again, repeating until there are no
more pairs of exposed carets with identical leaf numbers.  This
procedure produces the the unique {\em reduced} tree pair diagram
representing $g$.  When we write $g = (T,S)$, we are assuming that
this is the unique reduced tree pair diagram representing $g \in F$.
In this case, we refer to $T$ as the {\em negative} tree in the pair
and $S$ as the {\em positive} tree. This terminology is based on the
conversion of $(T,S)$ to the unique normal form of the element with
respect to the standard infinite generating set, and is described
explicitly in \cite{cfp}.

Let $T$ be a finite rooted binary tree with $n$ carets which we
number from $1$ through $n$ in infix order. We use the infix numbers
as names for the carets, and the statement $p<q$ for two carets $p$
and $q$ simply expresses the relationship between the infix numbers.
A caret is said to be a right (resp. left) caret if one of its sides
lies on the right (resp. left) side of $T$.  The root caret can be
considered either left or right.  All other carets are called
interior carets.

\begin{center}
\begin{figure}[ht]
  \includegraphics[scale=.75]{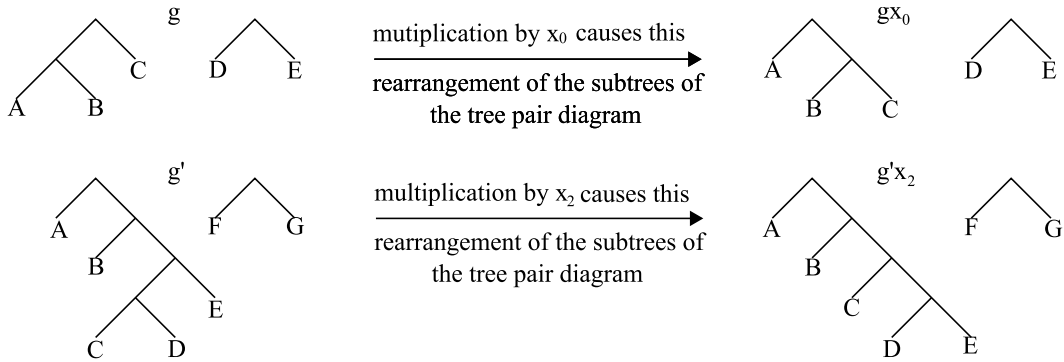}
  \caption{A depiction of the rearrangement of the subtrees of a tree pair diagram after
   multiplication by $x_0$ and $x_2$.  Capital letters represent (possibly empty) subtrees of the original tree
   pair diagram. Multiplication by $x_i^{\pm 1}$ causes an analogous rearrangement of the subtrees at level $i+1$.}
  \label{fig:multiplication}
\end{figure}
\end{center}

To multiply two elements $g = (T_1,T_2)$ and $h = (S_1,S_2)$ of $F$
we create unreduced representatives for the two elements, $g =
(T'_1,T'_2)$ and $h = (S'_1,S'_2)$ in which $S'_2 = T'_1$.  The
product $gh$ is then given by the (possibly unreduced) tree pair
diagram $(S'_1,T'_2)$.  In particular, if we take $h$ to be a
generator of the form $x_i^{\pm 1}$ we see that multiplication on
the right by $h$ causes a proscribed rearrangement of the subtrees
of $g = (T_1,T_2)$.  Note that it may be necessary to add carets to
the tree pair diagrams, creating unreduced representatives of these
elements, in order to preform this multiplication. The rearrangement
of the subtrees of $g$ under multiplication by $x_0^{\pm 1}$ and
$x_2^{\pm 1}$ is depicted in Figure \ref{fig:multiplication}.

Our formula for the word length of elements $g \in F$ with respect
to the generating set $X_n = \{x_0,x_1, \cdots ,x_n\}$ has two
components. The first we call $l_{\infty}(g)$, as it is the word
length of $g$ with respect to the standard infinite generating set
$\{x_i | i \geq 0\}$ for $F$.  This quantity is simply the number of
carets in the unique reduced tree pair diagram representing $g$
which are not right carets.  The second component in the word length
formula is twice what we term the {\em penalty weight} of the
element.  To make this precise, we begin by distinguishing a
particular type of caret in a single tree.

\begin{defn}[\cite{HST}, Definition 3.1]
Caret $p$ in a tree $T$ has {\em type N} if caret $p+1$ is an
interior caret which lies in the right subtree of $p$.
\end{defn}

We use this definition to describe certain carets in the tree pair
diagram for $g \in F$ which we call {\em penalty carets} as they
help determine the penalty contribution to the word length $l_n(g)$.
Let $g \in F$ have a reduced tree pair diagram $(T_-,T_+)$ in which
the carets are numbered in infix order.  By caret $p$ in $(T_-,T_+)$
we mean the pair of carets numbered $p$ in each tree.

\begin{defn}[\cite{HST}, Definition 3.2] Caret $p$ in a tree pair diagram $(T_-,T_+)$ is a {\it penalty caret} if either
\begin{enumerate}
\item $p$ has type $N$ in either $T_-$ or $T_+$, or
\item $p$ is a right caret in both $T_-$ and $T_+$ and caret $p$ is not the final caret in the tree pair
diagram.
\end{enumerate}
\end{defn}

To compute the penalty contribution to the word length for a given
$g=(T_-,T_+) \in F$ we use the following procedure.  Using a notion
of caret adjacency defined below, we take the two trees $T_-$ and
$T_+$ and construct a single tree $\PP$, called a {\em penalty
tree}, whose vertices correspond to a subset of the carets of $T_-$
and $T_+$, necessarily including the penalty carets.  This tree is
assigned a {\em weight} according to the arrangement of its
vertices. Minimizing this weight over all possible penalty trees
that can be constructed using the adjacencies between the carets of
$T_-$ and $T_+$ yields the penalty weight $p_n(g)$. We may now state
the word length formula precisely:

\begin{thm}[\cite{HST}, Theorem 3.3]\label{thm:length}
For every $g\in F$, the word length of $g$ with respect to the
generating set $X_n = \{x_0,x_1, \cdots ,x_n \}$ is given by the
formula
$$l_{X_n}(g)=l_n(g) = l_\infty(g)+2p_n(g)$$
where $l_{\infty}(g)$ is the number of carets in the reduced tree
pair diagram for $g$ which are not right carets, and $p_n(g)$ is the
penalty weight of $g$.
\end{thm}

\begin{center}
\begin{figure}[ht!]
  \includegraphics[scale=.75]{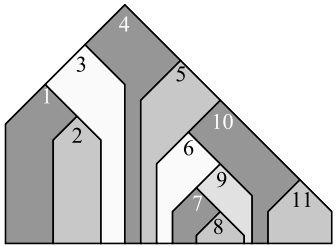}
  \caption{The shaded areas represent the carets of the tree, which are labeled in infix order.}
  \label{fig:adjacency}
\end{figure}
\end{center}

Constructing penalty trees for elements $g \in F$ requires a concept
of directed caret adjacency, which is an extension of the infix
order. To define the concept of adjacency between carets in a single
tree $T$, we view each caret as a space rather than an inverted v.
The point of intersection of the left and right edges of the caret
naturally splits the boundary of this space into a left and right
component. The space is bounded on the right (resp. left) by a {\em
generalized right (resp. left) edge}. The generalized right (resp.
left) edge may consist of actual left (resp. right) edges of other
carets in the tree, in addition to the actual right (resp. left)
edge of the caret itself. Let $p$ and $q$ denote carets in a tree
pair diagram $(T_-,T_+)$ and assume that $p < q$. We say that $p$ is
adjacent to $q$, written $p\prec q$, if there is a caret edge, in
either $T_-$ or $T_+$, which is both part of the generalized right
edge of caret $p$ and the generalized left edge of caret $q$. We
equivalently say that {\em traversing} the generalized left edge of
caret $q$ takes you to caret $p$ in at least one tree. It is always
true that carets $p$ and $p+1$ satisfy $p \prec p+1$. Although the
ordering of carets given by infix number is not symmetric but is
transitive, the notion of caret adjacency is neither symmetric nor
transitive.  Figure \ref{fig:adjacency} shows an example of a single
tree with the spaces corresponding to different carets shaded. In
this tree, in addition to the adjacency relationships $p \prec p+1$
for $1 \leq p \leq 10$, we also have $1 \prec 3$, $5 \prec 10$, $6
\prec 10$, $6 \prec 9$ and $7 \prec 9$.
We introduce a dummy caret denoted $v_0$ which is adjacent to all
left carets in both $T_-$ and $T_+$.  One can think of $v_0$ as
being the space to the left of the left side of each tree.  We now
construct a penalty tree $\PP$ corresponding to the pair of trees
$(T_-,T_+)$, which has this dummy caret $v_0$ as its root, according
to the following rules.

\begin{enumerate}
\item The vertices of $\PP$ are a subset of the carets in the tree pair
diagram, which we refer to by infix numbers: $0=v_0,1,2, \cdots ,k$,
always including $v_0$.
\item A directed edge may be drawn from vertex $p$ to vertex $q$ in
$\PP$ if $p \prec q$.
\item There is a vertex for every penalty caret in $(T_-,T_+)$.
\item Each leaf of $\PP$ corresponds to a penalty caret of
$(T_-,T_+)$.  The only exception to this is when $\PP$ consists only
of the root $v_0$ and no edges.
\end{enumerate}
The penalty tree $\PP$ is oriented in the sense that there is a
unique path from $v_0$ to every vertex $p \in \PP$, and if this path
passes through vertices $v_0, p_1, p_2, \ldots, p_{i}=p$ then we
must have $v_0 \prec p_1 \prec \cdots \prec p_i=p$.  Two vertices
$p, q$ in the tree are comparable if there is either a path $p=w_1,
w_2, \ldots, w_{i+1}=q$ or $q=w_1, w_2, \ldots, w_{i+1}=p$ with $w_j
\prec w_{j+1}, \forall j=1, \ldots, i+1$, and in this case we say
$d_{\PP}(p,q)=i$.

The penalty weight of a penalty tree is bounded above by the number of vertices on the tree, but not all vertices on the tree contribute to the weight. More precisely,
we define:
\begin{defn}[\cite{HST}, Definition 3.4]The {\em n-penalty weight} $p_n(\PP)$ of a penalty tree $\PP$
associated to $g = (T_-,T_+) \in F$ is the number of vertices $v_i
\in \PP$ such that $d_{\PP}(v_0, v_i) \geq 2$ and there exists a
leaf $l_i$ in $\PP$ with $d_{\PP}(v_i, l_i) \geq n-1$. These
vertices are called the {\em weighted} carets.
\end{defn}

To compute the penalty contribution $p_n(g)$ to the word length
$l_n(g)$ for $g \in F$, we must minimize the penalty weight over all
penalty trees associated to $g$.
\begin{defn}[\cite{HST},Definition 3.5] For an element $g \in F$, define the penalty weight of the element $G \in F$,
denoted $p_n(g)$ by $$p_n(g)= min \{p_n(\PP)| \PP \mbox{ is a
penalty tree for } g=(T_-,T_+) \}$$
\end{defn}

We have now defined both components of the word length formula given
in Theorem \ref{thm:length}.

\section{Coarse isometries and convexity}
Recall that a map $f$ between two metric spaces $G$ and $H$ is a
\emph{quasi-isometry} if there are positive constants K and C so
that for every pair of points $g_1, g_2 \in G$,
 $$\frac{1}{K}d_G(g_1, g_2)-C \leq d_H(f(g_1),f(g_2)) \leq Kd_G(g_1,g_2)+C.$$
If the constant $K$ can be chosen to be $1$, we call $f$ a
\emph{C-coarse isometry}. Given a group $G$ and a finite generating
set $X$, G can be regarded as a metric space using the wordlength
metric, namely, $d_X(g,h)=min\{n|gh^{-1}=\alpha_1 \alpha_2 \cdots
\alpha_n, \alpha_i^{\pm1} \in X \}$. We denote $G$, viewed as a
metric space in this way, by $(G,X)$. Equivalently, one can view the
Cayley graph $\Gamma(G,X)$ as a metric space by declaring each edge
to have length 1. Recall that for any finitely generated group $G$
with finite generating sets $X$ and $Y$, the identity map between
$(G,X)$ as $(G,Y)$ is a quasi-isometry. In general, it is unknown to
what extent quasi-isometries preserve convexity properties, but in
the special case of a coarse-isometry, we obtain the following:

\begin{theorem}
\label{thm:notmac} Let $f:(G,X_G) \rightarrow (H,X_H)$ be a
$C$-coarse-isometry. If $(G,X_G)$ is not minimally almost convex,
then $(H,X_H)$ is not almost convex.
\end{theorem}

\begin{proof}
Let $g$ be any coarse inverse for $f$, which is easily seen to be a
coarse isometry as well.  Without loss of generality, we may assume
that $g$ is also a $C$-coarse isometry.

Suppose that $(H,X_H)$ is almost convex.  Then for each $n\geq 2$,
there is an almost convexity constant $K(n)$.  Fix $M > 2C+1$, and
let $K = K(2M + C)$.  Let $n > K + M + C$.

Since $(G,X_G)$ is not minimally almost convex, we can find $x,y \in
B(n) \subset \Gamma(G,X_G)$ with $d_G(x,y) = 2$ so that the shortest
path from $x$ to $y$ which remains in $B(n)$ has length $2n$. Since
we can always construct a path of this length passing through the
identity, let $\gamma$ be such a path containing the identity.

Consider the closed loop $\eta$ obtained by concatenating $\gamma$
with the path of length two between $x$ and $y$.  Let $z$ denote the
point in $B(n+1)$ at distance one from $x$ and $y$.  Choose $a$ and
$b$ on $\gamma$, with $a$ on the subpath of $\gamma$ from $x$ to the identity, and $b$ between $y$ and the identity, so that $d_G(a,Id) = d_G(b,Id)$ and $d_G(a,z) =
d_G(b,z) = M$.  Let $\eta_1$ be the subpath of $\gamma$ containing
$a, \ b$ and the identity, and $\eta_2$ is the remaining subpath of
$\eta$.

Consider $f(a)$ and $f(b)$, elements of the Cayley graph
$\Gamma(H,X_H)$.  We know that $d_H(f(a),f(b)) \leq 2M + C$. Since
we are assuming that $(H,X_H)$ is almost convex, there must be a
path $\xi$ from $f(a)$ to $f(b)$ whose length is at most $K$, and
which remains in the ball $B(D)$, where $D$ is defined by $D=\max
\{d_H(f(a),id), \ d_H(f(b),id)\} \leq d_G(a,id) + C$.

Consider the image of $\xi$ under $g$, the coarse inverse to $f$.
Since $length(\eta_1)=2n-2M+2 > 2(K+C+M)-2M+2 > 2K+2C$ and
$length(g(\xi)) \leq K+C$, we see that $length(g(\xi)) <
length(\eta_1)$.  We now show that this path stays in $B(n)$,
contradicting the fact that any path from $x$ to $y$ in $B(n)$ has
length $2n$.

The maximum distance of any point on $\xi$ from the identity in $H$
is $D$.  Thus the maximum distance of any point on $g(\xi)$ from the
identity of $G$ is $D+C \leq d_G(a,id) + 2C=n-M+1+2C$.  Since $M
>2C+1$, it follows that $g(\xi) \subset B_G(n)$.

By concatenating the portion of $\eta_2$ from $x$ to $a$, $g(\xi)$,
and the portion of $\eta_2$ from $b$ to $y$, we obtain a path from
$x$ to $y$ which remains inside of $B(n)$ and has length less than
$2n$, a contradiction since $(G,X_G)$ is not minimally almost
convex. \end{proof}

\subsection{Application to Thompson's group $F$}

In \cite{CT1} it is shown that Thompson's group $F$ is not almost
convex with respect to the standard finite generating set
$\{x_0,x_1\}$. A natural question is whether $F$ is not almost
convex with respect to any finite generating set.  We use Theorem
\ref{thm:notmac} to extend this result to finite generating sets for
$F$ of the form $\{x_0,x_n\}$.

Belk and Bux in \cite{BBu} show that Thompson's group $F$ is not
minimally almost convex with respect to the generating set
$\{x_0,x_1\}$.  It is easy to see that the the word metrics in $(F,
\{x_0,x_1\})$ and$(F,\{x_0,x_n\})$ differ by the additive constant
$2(n-1)$, and thus the quasi-isometry between these two
presentations for $F$ is a coarse isometry.

Combining these results with Theorem \ref{thm:notmac}, we obtain the
following corollary, which is a special case of Theorem
\ref{thm:notacgeneral} below.

\begin{corollary} \label{cor:F} Thompson's group $F$ is not almost convex with
respect to any generating set of the form $\{x_0,x_n\}$.
\end{corollary}

\section{Convexity results}
The main goal of this section is to show that $F$ is not almost
convex with respect to any generating set which is a subset of the
standard infinite generating set; we note that in order for a subset
of the standard infinite generating set to generate $F$, it must
contain $x_0$. We extend the result of \cite{BBu} which proves that
$(F,\{x_0,x_1\})$ is not minimally almost convex first to
consecutive generating sets for $F$, and then to generating sets
which contain $x_0$ and $x_1$ and are subsets of the standard
infinite generating set.  To obtain our ultimate result, that $F$ is
not almost convex with respect to any generating set which is a
subset of the infinite generating set for $F$ and contains $x_0$, we
again discuss coarse isometries between different presentations for
$F$.

We begin with the following:
\begin{theorem}\label{thm:notMACconsec}
Let $X_n = \{x_0,x_1, \cdots ,x_n\}$ be a consecutive generating set
for $F$ with $n \geq 2$.  Then $F$ is not minimally almost convex with respect to
$X_n$.
\end{theorem}

\begin{proof}
We prove this by providing, for any $k>0$, a pair of group elements
$g=g_k$ and $h=h_k$ satisfying $l_n(g)=l_n(h)= 2k+2$ and
$l_n(h^{-1}g)=2$, for which any path $\gamma$ from $g$ to $h$ that
lies entirely within the ball of radius $2k+2$ must have length at
least $4k+4$.

Let $g = g_n =x_1^{k+1}x_{k+n+1}x_0^{-k}=x_nx_1^{k+1}x_0^{-k}$ and
$h = h_n = gx_0^{-1}x_n^{-1}=x_1^{k+1}x_0^{-(k+1)}$.  The tree pair
diagrams for these elements are given in Figure \ref{fig:g}. In the
tree pair diagrams for $g$ and $h$, we observe that
$l_\infty(g)=l_{\infty}(h)=2k+2$.  From Theorem \ref{thm:length} we
see that $l_n(a) \geq l_{\infty}(a)$ for all $a \in F$, and since we
have provided words above of length $2k+2$ for both $g$ and $h$, it
follows that $l_n(g)=l_n(h)= 2k+2$.

\begin{center}
\begin{figure}[ht]
  \includegraphics[scale=1]{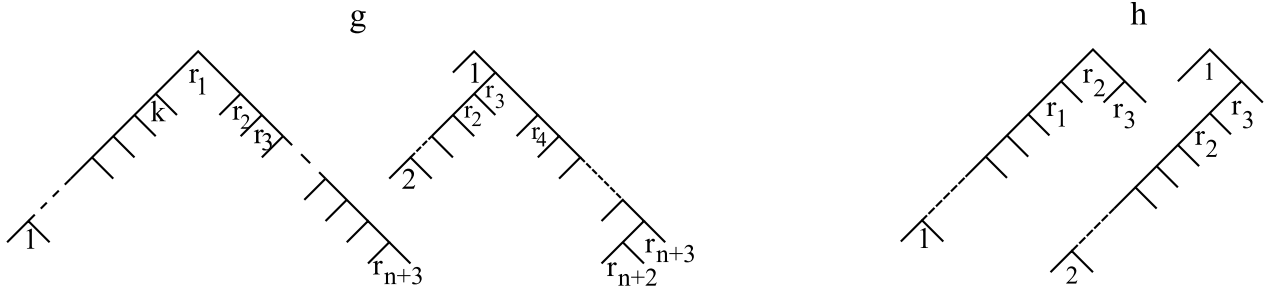}
  \caption{The tree pair diagrams representing the elements $g$ and $h$ used in the proof of
  Theorem \ref{thm:notMACconsec}.}
  \label{fig:g}
\end{figure}
\end{center}

Suppose there is a path $\gamma$ from $g$ to $h$ which lies within
the ball of radius $2k+2$.  We note that the only generator $x \in
X_n$ so that the word length of $gx$ is less than the word length of
$g$ is $x = x_0$.  Thus the first vertex along $\gamma$ after $g$ is
$gx_0$.  In the negative tree for the tree pair representing $gx_0$,
the caret $r_{n+2}$ is a right caret at level $n+3$, whereas in the tree pair
diagram for $g$ it is a right caret at level $n+2$.  Our argument relies on noting
the level of this caret at successive vertices along the path
$\gamma$.

In order for the path $\gamma$ to terminate at $h$, there is a point
at which the pair of carets numbered $r_{n+2}$ in each tree must be
removed as part of a reduction along $\gamma$.  This requires caret
$r_{n+2}$ from $T_-$ to be an interior caret at the point of
reduction. Given the effect of
multiplication by each generator on the tree pair diagram as
described in Section \ref{sec:computinginF}, we observe that the generators in $X_n$ cannot move any right caret off the right side of the tree unless it is at level 1 through $n+1$. Hence, we conclude that there is a
smallest nontrivial prefix $\gamma_0$ of $\gamma$ so that in $g
\gamma_0 = f$ the caret $r_{n+2}$ in the negative tree for $f$ is a right caret at
level $n+1$.

Let $(S_-,S_+)$ be the tree pair diagram for $f = g \gamma_0$ which
is constructed from the tree pair diagram $(T_-,T_+)$ for $g$ by
altering these trees according to multiplication by each generator
of $\gamma_0$, but without performing any possible reductions.
During this process, the carets in $T_+$ remain unchanged, though
additional carets may be added to $T_+$ to form $S_+$. Hence, $S_+$
contains $T_+$ as a subtree, and the tree pair diagram $(S_-,S_+)$
may be unreduced.

We first show that the tree pair diagram $(S_-,S_+)$ constructed in
this way must be unreduced, and that when the reduction is
accomplished, some of the original carets from
$T_+$ will be removed from $S_+$.  If this was not the case, then in
$S_-$ there would be at least $k+1$ carets with smaller infix
numbers than $r_{n+1}$ which were not right carets, and thus counted
towards $l_{\infty}(f)$.  Additionally, in $S_+$ there would also be $k+1$
interior carets with infix numbers less than $r_{n+2}$, and caret
$r_{n+2}$ itself is also an interior caret.  This implies that
$l_{\infty}(f) \geq 2k+3$, contradicting the fact that $f \in
B(2k+2)$.  Thus there must be some reduction of the carets of
$T_+$, viewed as a subtree of $S_+$, in order to obtain the reduced tree pair diagram for $g \gamma_0 =f$.

We now consider which carets of $T_+$, viewed as a subtree of $S_+$
might be reduced; in order for a caret to be reduced after
multiplication by a particular generator, it must be exposed, that
is, both leaves have valence one. The only exposed carets of $T_+$
itself are carets $2$ and $r_{n+2}$.  Since caret $r_{n+2}$ is a
right caret in $S_-$, and not the final right caret, it is not
exposed in $S_-$. Therefore, it must be that in reducing
$(S_-,S_+)$, the original caret $2$ from the infix ordering on $T_+$
must cancel. We claim that in $S_-$, caret $2$ must be a child of
caret $1$. If, in forming $S_+$, no carets were added to either leaf
of caret $2$, then caret $2$ is exposed in $S_+$, and hence it is
exposed in $S_-$, which implies that caret $2$ is a child of caret
$1$ in $S_-$. If, on the other hand, carets were added to the leaves
of caret $2$ in forming $S_+$, then they must all cancel in
$(S_-,S_+)$ before caret $2$ does. But this means that in $S_-$,
these added carets must also hang from the leaves of caret $2$, and
once again, caret $2$ is a child of caret $1$ in $S_-$.

The fact that caret $2$ is a child of caret $1$ in $S_-$ provides a lower bound on $l_n(h^{-1}f)$ as
follows.  To form the tree pair diagram for $h^{-1}f$, consider the
unreduced tree pair diagram $(S_-,S_+)$.  If $h = (H_-,H_+)$, to
form this product we consider these trees in the order $S_- \ S_+ \
H_+  \ H_-$, and add carets to each pair to ensure that the middle
trees are identical.  Thus we must at least add the string of right carets
$r_4, \ldots, r_{n+1}, r_{n+3}$, with caret $r_{n+2}$ the left child
of $r_{n+3}$, from $S_+$ to both trees in the diagram $(H_+,H_-)$ in
order to perform this multiplication. Since in $S_-$, caret $2$ is a child of caret $1$, but in $H_-$ caret $1$ is a child of caret $2$, caret $1$ cannot reduce in the product $h^{-1}f$. Hence, because of their configuration in $H_-$, the entire string of carets $1,2, \cdots,k ,r_1$ do not
reduce in the product $h^{-1}f$. Also, as we remarked above, caret
$r_{n+2}$ is not removed through reduction in this product.  Hence
we obtain the following lower bound on the word length of $h^{-1}f$:
$l_n(h^{-1}f)\geq l_{\infty}(h^{-1}f) \geq 2(k+1)+1=2k+3.$

Let $\gamma_1$ be the subpath of $\gamma$ from $f = g\gamma_0$ to
$h$.  Since $l_n(h^{-1}f) \geq 2k+3$, it follows that $|\gamma_1|
\geq 2k+3$. But traversing $\gamma_0$ in reverse, followed by
$x_0^{-1}$ and then $x_n^{-1}$ yields another path from $f$ to $h$,
so similarly $|\gamma_0|+2 \geq 2k+3$, and hence $|\gamma_0| \geq
2k+1$. This implies that $|\gamma|=|\gamma_0|+|\gamma_1| \geq 4k+4$.
\end{proof}

In the proof above, both $g$ and $h$ are be represented by words of
length $2k+2$ involving only the generators $x_0^{\pm 1}, x_1^{\pm
1}$, and $x_n^{\pm 1}$, namely, $g=x_nx_1^{k+1}x_0^{-k}$ and
$h=x_1^{k+1}x_0^{-(k+1)}$. Hence, the above result can be extended
to any generating set for $F$ which is a finite subset of the
standard infinite generating set containing $x_0$ and $x_1$.
\begin{theorem}\label{thm:notMACwithx_1}
Let $X=\{x_0,x_1, x_{i_1},x_{i_2}, \cdots ,x_{i_j}\}$, where $1< i_1 < \cdots < i_j$, be a
generating set for $F$.  Then $F$ is not minimally almost convex
with respect to $X$.
\end{theorem}

\begin{proof}
The identity map on $G$ is a quasi-isometry between the metric
spaces $(G,X)$ and $(G,X_{i_j})$, where $X_{i_j} = \{x_0,x_1,x_2,x_3
\cdots ,x_{i_j}\}$. Since $X \subset X_{i_j}$, we remark that
$d_{X_{i_j}}(a,b) \leq d_X(a,b)$ for any $a,b \in F$. In particular,
$d_{X_{i_j}}(a,Id)) \leq d_X(a,Id)$ for any $a \in F$.

Assume that $(F,X)$ is minimally almost convex.  It is proven in
Theorem \ref{thm:notMACconsec} that $(F, X_{i_j})$ is not minimally
almost convex.  Let $h=h_k=x_1^{k+1}x_0^{-(k+1)}$ and $g=g_k=
x_1^{k+1}x_{k+i_j+1}x_0^{-k}$  be the group elements used in the
proof of Theorem \ref{thm:notMACconsec}.  It is clear that
$2k+2=d_{X_{i_j}}(h,id) = d_X(h,id)$ and $2k+2=d_{X_{i_j}}(g,id) =
d_X(g,id)$; if there was a shorter expression for either $g$ or $h$
with respect to $X$, then there would be one with respect to
$X_{i_j}$ as well. In addition, it is clear that since
$g^{-1}h=x_{i_j+1}^{-1}x_0^{-1} = x_0^{-1}x_{i_j}^{-1}$, we have
$d_{X_{i_j}}(g,h)=d_X(g,h) = 2$.

Since $(F,X)$ is assumed to be minimally almost convex, there is a
path $\gamma$ of length at most $4k+3$ connecting $g$ and $h$ which
lies within the ball of radius $2k+2$ relative to $X$. Since each
group element $a$ along this path satisfies $d_{X_{i_j}}(a,id)) \leq
d_X(a,id) \leq 2k+2$, this contradicts the assumption that
$(F,X_{i_j})$ is not minimally almost convex.  Thus we conclude that
$(F,X)$ cannot be minimally almost convex.
\end{proof}

To extend the result of Theorem 6.1 of \cite{HST} to arbitrary
finite subsets of the infinite generating set containing $x_0$, we
show first that word length with respect to one of these arbitrary
generating sets differs from word length with respect to some
generating set containing $x_1$ only by an additive constant.

\begin{lemma}
Let $X=\{x_0,x_{i_1},x_{i_2}, \cdots ,x_{i_j}\}$ be a generating
set for $F$, and form a new generating set $Y =\{x_0, x_1, x_{i_2-i_1+1}, x_{i_3-i_1+1}, \ldots x_{i_j-i_1+1}\}$.
Then $(F,X)$ and $(F,Y)$ are coarsely isometric.
\end{lemma}
\begin{proof} Let $g \in F$, and suppose $g= \alpha_1 \alpha_2 \cdots \alpha_m$, where
$\alpha_k^{\pm1} \in Y$.  Then $$g=x_0^{i_1-1} \left( x_0^{1-i_1} g
x_0^{i_1-1}\right) x_0^{1-i+1} = x_0^{i_1-1} \bar{\alpha_1}
\bar{\alpha_2} \cdots \bar{\alpha_m} x_0^{i_1-1},$$ where
$\bar{\alpha_k}=x_0^{1-i_1}\alpha_k x_0^{i_1-1}$. Now in the cases
where $\alpha_k=x_0^{\pm1}$, we have $\bar{\alpha_k}=\alpha_k$, and
in the cases where $\alpha_k=x_l^{\pm1}$ with $l \geq 1$, then
$\bar{\alpha_k}= x_{l+i_1-1}^{\pm 1}\in X$. Hence $l_{X}(g) \leq
l_{Y}(g)+2(i_1-1)$. Similarly, one sees that $l_{Y}(g) \leq
l_X(g)+2(i_1-1)$. Hence, $l_{X}(g) -2(i_1-1) \leq l_{Y}(g) \leq
l_{X}(g)+2(i_1-1)$ and $l_{Y}(g) -2(i_1-1) \leq l_{X}(g) \leq
l_{Y}(g)+2(i_1-1)$.
\end{proof}
Finally, we apply Theorem \ref{thm:notmac} to $F$ with the two generating sets $X$ and $Y$ of the preceding theorem to obtain:
\begin{theorem}\label{thm:notacgeneral}
Let $X$ be any subset of the standard infinite generating set for
$F$ which includes $x_0$.  Then $F$ is not almost convex with
respect to $X$.
\end{theorem}

\end{document}